\documentclass[12pt]{article}
\usepackage{amsmath,amssymb}

\def\R{{\mathbb R}}

\def\bf{{\bfseries}}

\def\A{\mathbb A}

\def\st{\vphantom{{(\over)}}}
\def\cc#1={\cite{#1}}
\def\fff#1={\eqno (#1)}
\def\text#1{\hbox{#1}}

\def\deff{\buildrel\rm def \over =}
\def\fr#1#2{{#1 \over #2}}
\def\sfr#1#2{{\textstyle{ #1 \over #2 }}}
\def\ttt{\hbox{\ {.}\,{.}\,{.}\ }}

\def\nad#1#2{\buildrel{#1} \over{#2}\!\!\strut}
\def\pod#1#2{\mathrel{\mathop{#2}\limits_{#1}}\strut}

\def\eq{\equiv}
\def\neq{\not\equiv}
\def\nr{\ne}
\def\mr{\leq}

\def\ds{\displaystyle}

\def\ll{${\cal l}$}
\def\an{$M_n$}
\def\ann{$\bar M_n$}
\def\kn{$K_n$}
\def\knn{$\bar K_n$}
\def\vn{$V_n$}
\def\vnn{$\bar V_n$}

\def\dd{\partial}
\def\Rank{\,\hbox{\rm Rank}\,}

\def\a{\alpha}
\def\vph{\vphantom{\alpha}}
\def\b{\beta}
\def\f{\varphi}
\def\Gf{\varPhi}
\def\c{\psi}
\def\Gc{\Psi}
\def\g{\gamma}
\def\Gg{\Gamma}
\def\d{\delta}
\def\e{\varepsilon}
\def\l{\lambda}
\def\Gl{\Lambda}
\def\go{\omega}
\def\Go{\Omega}
\def\s{\sigma}
\def\Th{\Theta}
\def\ep{\varepsilon}
\def\phi{\varphi}
\def\r{\varrho}

\newtheorem{theorem}{Theorem}

\newtheorem{coro}{Corollary}

\date{}
\begin{document}

\title{FUNDAMENTAL PDE'S\\
OF THE CANONICAL ALMOST  \\ GEODESIC MAPPINGS OF TYPE ${\tilde\pi}_1$}

\author{V. E.~Berezovski, J.~Mike\v s, A.~Van\v zurov\'a
\vphantom{
\thanks{Supported by the Council of Czech Government MSM 6198959214.}
\thanks{  {\it Author's addresses:\/}
Dept.~of~Mathematics, University of Uman, Institutskaya~1, Uman, Ukraine;
Dept.~of Algebra and Geometry, Fac.~Sci., Palack\'y Univ., T\v r.~17. listopadu 12,
        770 00 Olomouc, Czech Rep.,
e-mail: berez.volod@rambler.ru; josef.mikes@upol.cz; vanzurov@inf.upol.cz.
            }
        }
            }

\maketitle

\begin{abstract}
For modelling of various physical processes, geodesic lines and almost geodesic curves serve as a useful tool. Trasformations or mappings between spaces (endowed with the metric or connection) which preserve such curves play an important role in physics, particularly in mechanics, and in geometry as well.
Our aim is to continue investigations concerning existence
of almost geodesic mappings of manifolds with linear (affine) connection, particularly
of the so-called ${\tilde \pi}_1$ mappings, i.e.~canonical almost geodesic mappings
of type $\pi_1$ according to Sinyukov.
First we give necessary and sufficient conditions for existence of ${\tilde \pi}_1$ mappings
of a manifold endowed with a linear connection onto
pseudo-Riemannian manifolds. The conditions take the form
of a closed system of PDE's of first order of Cauchy type.
Further we deduce necessary and sufficient conditions for existence of
${\tilde \pi}_1$ mappings onto generalized Ricci-symmetric spaces.
Our results are generalizations of some previous theorems obtained by N.S.~Sinyukov.
\end{abstract}

\medskip


%
%
\def\st{\vphantom{{(\over)}}}
\def\cc#1={\cite{#1}}
\def\fff#1={\eqno (#1)}
\def\text#1{\hbox{#1}}
\def\deff{\buildrel\rm def \over =}
\def\fr#1#2{{#1 \over #2}}
\def\sfr#1#2{{\textstyle{ #1 \over #2 }}}
\def\ttt{\hbox{\ {.}\,{.}\,{.}\ }}
\def\nad#1#2{\buildrel{#1} \over{#2}\!\!\strut}
\def\pod#1#2{\mathrel{\mathop{#2}\limits_{#1}}\strut}
\def\eq{\equiv}
\def\neq{\not\equiv}
\def\nr{\ne}
\def\mr{\leq}
\def\ds{\displaystyle}
\def\ll{${\cal l}$}
\def\an{$M_n$}
\def\ann{$\bar M_n$}
\def\kn{$K_n$}
\def\knn{$\bar K_n$}
\def\vn{$V_n$}
\def\vnn{$\bar V_n$}
\def\dd{\partial}
\def\Rank{\,\hbox{\rm Rank}\,}
\def\a{\alpha}
\def\vph{\vphantom{\alpha}}
\def\b{\beta}
\def\f{\varphi}
\def\Gf{\varPhi}
\def\c{\psi}
\def\Gc{\Psi}
\def\g{\gamma}
\def\Gg{\Gamma}
\def\d{\delta}
\def\e{\varepsilon}
\def\l{\lambda}
\def\Gl{\Lambda}
\def\go{\omega}
\def\Go{\Omega}
\def\s{\sigma}
\def\Th{\Theta}
\def\ep{\varepsilon}
\def\phi{\varphi}
\def\r{\varrho}
%

\section{Introduction}

\strut\indent
Geodesic and almost geodesic lines serve as a useful tool for modelling of various physical processes, and mappings between spaces
(endowed with the metric or connecion) and trasformations which preserve such curves, play an important role in geometry as well as in physics, particularly in mechanics, optics and the theory of relativity, \cite{Pe}, \cite{Pe1}, \cite{Pe2}.

Many geometric problems connected with the topic of differential ge\-o\-met\-ry are solved by means of
differential equations, particularly, the problems are often answered by solving systems of
partial differential equations (PDE's) for components of some geometric objects (e.g.~tensors),
\cite{Ei,Fe,MiKV,MVH,Pe,Si79}.
We intend to study here the existence problem of canonical almost geodesic mappings, and as we shall see,
our main tool will be to construct and solve a suitable system of PDE's of Cauchy type
that controlles the situation. One of characteristic properties of a system of PDE's of Cauchy type is
that the solution of such a system depends on a finite number of real (or complex) parameters. Moreover,
solutions of such systems can be effectively enumerated, eventually some approximation can be found.

Unless otherwise specified, all objects under consideration are supposed
to be differentiable of a sufficiently high class (mostly, differentiability
of the class $C^3$ is sufficient).

Let $A_n=(M,\nabla)$ be an $n$-dimensional ($C^{k}$, $C^{\infty}$ or $C^{\omega}$)
manifold endowed with a linear connection $\nabla$.
Let $c:I\to M$, $t\mapsto {c}(t)$ defined on an open interval
$I\subset\R$ be a ($C^k$, or smooth) curve on $M$ satisfying the
regularity condition
\begin{displaymath}
c'(t)=
dc(t)/dt\ne 0 \mbox{\rm\quad for all\ \ }t\in I.
\end{displaymath}
Denote by $\xi$ the corresponding ($C^{k-1}$, or smooth)
tangent vector field along~$c$ (``velocity field"),
$ 
\xi(t)=\left (c(t),c'(t)\right), \ t\in I,
$ 
and let
\begin{equation}
{\xi}_1=\nabla(\xi;\xi)={\nabla}_{\xi} \xi,\qquad {\xi}_2=
{\nabla}^2(\xi;\xi,\xi)= {\nabla}_{\xi} {\xi}_1.
\end{equation}
Geodesics $c(s)$, parametrized by canonical affine parameter (given up
to the affine transformations $s\mapsto as+b$),
are characterized by ${\nabla}_{\xi} {\xi}=0$ while
unparametrized geodesic curves (i.e.~arbitrarily parametrized,
called also \emph{pregeodesics} in the literature)
can be characterized
by the formula ${\nabla}_{\xi} {\xi}={\lambda}{\xi}$
where $\lambda(t):I\to\R$ is a real function.

Let $D=\mbox{\rm span\,}(X_1,X_2)$ (i.e.\ the vector fields
$X_1$, $X_2$ along $c$ form a basis of $D$). Recall that $D$ is parallel
(along $c$) if and only if the covariant derivatives along $c$
of basis vector fields belong to the distribution
(the property is independent of reparametrization of the curve)
\cite{Si63,Si79,Si82}.

As a generalization of (an unparametrized) geodesic,
let us introduce an \emph{almost geodesic curve} as
a curve $c$ satisfying: there exists a two-dimensional
(differentiable) distribution $D$
parallel along $c$ (relative to $\nabla$)
such that for any tangent vector of $c$, its parallel
translation along $c$ (to any other point)
belongs to the distribution $D$.
Equivalently, $c$ is almost geodesic if and only if
there exist vector fields $X_1$, $X_2$ parallel along $c$ (i.e.\
satisfying $\nabla_{\xi} X_i=a^jX_j$ for some differentiable
functions $a^j_i(t):I\to\R$)
and differentiable real functions $b^i(t)$, $t\in  I$ along $c$,
such that $\xi=b^1X_1+b^2X_2$ holds.
For almost geodesic curves, the vector fields ${\xi}_1$ and ${\xi}_2$
belong to the corresponding distribution $D$.
If the vector fields ${\xi}$ and ${\xi}_1$ are independent at any point
(and hence the (local) curve $c$ is not a geodesic one),
we can write $D=\mbox{\rm span\,}({\xi},{\xi}_1)$.
So we get another equivalent characterization:
{\it a curve is almost geodesic if and only
if ${\xi}_2\in \mbox{\rm span\,}(\xi,{\xi}_1)$.\/}

\section{Almost geodesic mappings}

Geodesic mappings of mani\-folds with linear connection are ($C^k$)-diffeo\-mor\-phisms
characterized by the property that all geodesics are send onto
(unparametrized in general) geodesic curves. The classification of
geodesic mappings is more or less known.
Recall that even for Riemannian spaces, there is a lack of a nice simple
criterion for decision when a given Riemannian space admits non-trivial
geodesic mappings.

Let $A_n=(M,\nabla)$, ${\bar A}_n=(\bar M,\bar{\nabla})$ be
$n$-dimensional manifolds ($n>2$) each endowed with a torsion-free
linear connection.

We may ask which ($C^k$-)diffeomorphisms of manifolds
send almost geodesic curves onto almost geodesic again. The answer is: such mappings reduce to
geodesic ones, since there are ``too many" almost geodesic curves.
It appears that the following definition is more acceptable.

We say that a ($C^k$-)diffeomorphism $f$: $M\to \bar M$ is \emph{almost geodesic}
if any geodesic curve of $(M,\nabla)$
is mapped under $f$ onto an almost geodesic curve in
$(\bar M,\bar{\nabla})$.

This concept of an almost geodesic mapping
was introduced by N.S.~Si\-nyu\-kov
\cc Si63=,
and before by V.M.~Chernyshenko
\cc Che=,
from a rather different point of view.
The theory of almost geodesic mappings was treated in
\cc Si63, Si79, Si82=.

Due to the fact that $f$ is a diffeomorphism we can accept the useful
convention that both linear connections $\nabla$ and $\bar{\nabla}$ are
in fact defined on the same underlying manifold $M$,
so that we can consider their
difference tensor field of type $(1,2)$, $P=\bar\nabla-\nabla$,
called sometimes a \emph{deformation tensor} of the given connections
under $f$ \cc Si79=,
given by $\bar\nabla(X,Y)=\nabla(X,Y)+P(X,Y)$ for $X,Y\in{\cal X}(M)$.
Since the connections are symmetric, $P$ is also symmetric in $X,Y$.
Of course, we identify objects on $M$ with their corresponding objects on
$\bar M$: a curve $c$  on $M$ identifies with its image
$\bar c=f\circ c$, its tangent vector field $\xi(t)$ with the corresponding
vector field $\bar\xi(t)=Tf(\xi(t))$ etc.

Besides the deformation tensor, we will use the tensor field of type $(1,3)$,
denoted by the same symbol $P$, introduced by
\begin{displaymath}
P(X,Y,Z)= \!\!\! \mathop{\sum}_{CS(X,Y,Z)} \!\!\!
{\nabla}_ZP(X,Y)+P(P(X,Y),Z),\quad X,Y,Z\in{\cal X}(M),
\end{displaymath}
where
${\sum}_{CS(\, ,\, ,\, )}$
means the cyclic sum on
arguments in brackets
(i.e.~symmetrization without coefficients).

Almost geodesic diffeomorphisms $f\colon (M,{\nabla})\to (M,\bar{\nabla})$
are characterized by the following condition on the type $(1,3)$ tensor
$P$:
\begin{displaymath}
P(X_1,X_2,X_3)\wedge P(X_4,X_5)\wedge X_6=0, \quad
X_i\in{\cal X}(M),\ \ i=1,\dots ,6;
\end{displaymath}
$X\wedge Y$ means the decomposable bivector, the exterior product of $X$ and~$Y$.

N.S.~Sinyukov \cc Si63, Si79, Si82= distinguished three kinds of almost geodesic mappings,
namely
${\pi}_1$, ${\pi}_2$, and ${\pi}_3$, characterized, respectively,
by the conditions for the deformation tensor:
\begin{displaymath}
{\pi}_1{:}\ \nabla_X P(X,X)+P(P(X,X),X)
=a(X,X){\cdot}X+b(X){\cdot} P(X,X),\
X\in{\cal X}(M),\label{2.4}
\end{displaymath}
where $a\in S^2(M)$ is a symmetric tensor field  of type $(0,2)$ and $b$
is a $1$-form;
\begin{displaymath}
{\pi}_2{:}\ P(X,X)=\psi(X)\cdot X+\varphi(X)\cdot F(X),\quad
X\in{\cal X}(M),
\end{displaymath}
where $\psi$ and $\phi$ are $1$-forms, and $F$ is a type $(1,1)$
tensor field satisfying
\begin{displaymath}
\nabla_X F(X)+\phi(X)\cdot F(F(X))=\mu(X)\cdot X+\varrho(X)\cdot F(X),\quad
X\in{\cal X}(M)
\end{displaymath}
for some $1$-forms $\mu$, $\varrho$;
\begin{displaymath}
{\pi}_3{:}\ P(X,X)=\psi(X)\cdot X+a(X,X)\cdot Z,\quad
X\in{\cal X}(M)
\end{displaymath}
where $\psi$ is a $1$-form, $a\in S^2(M)$ is a symmetric bilinear form and
$Z\in{\cal X}(M)$ is a vector field satisfying
\begin{displaymath}
\nabla_X Z=h\cdot X+\theta(X)\cdot Z
\end{displaymath}
for some scalar function $h$: $M\to \R$ and some $1$-form $\theta$.
Remark that the above classes are not disjoint.

\section{Canonical almost geodesic mappings ${\tilde{\pi}}_1$}

We are interested here in a particular subclass of ${\pi}_1$-mappings,
the so-called $\tilde{\pi}_1$-{\it mappings},
or \emph{canonical} almost geodesic mappings, distinguished by the
condition $b=0$.
That is, ${\tilde{\pi}}_1$-mappings are just morphisms satisfying
\begin{displaymath}
\nabla_X P(X,X)+P(P(X,X),X)=a(X,X)\cdot X, \quad a\in S^2(M),\
X\in{\cal X}(M).
\end{displaymath}
In local coordinates, the condition reads
\begin{equation}
P^h_{(ij,k)}=a_{(ij}{\delta}^h_{k)}-P^h_{\alpha (i} P^{\alpha}_{jk)}.
\label{2.9}
\end{equation}
Here and in what follows, 
the comma ``\,,\," denotes covariant derivative with respect to $\nabla$,
$\d^h_i$ is the Kronecker delta, the round bracket denote the cyclic sum on indices
involved.

Any geodesic mapping is a ${\pi}_1$-mapping (the characterizing
condition can be checked), and any ${\pi}_1$-mapping can be
written as a composition of a geodesic mapping followed by a
${\tilde{\pi}}_1$-mapping.
So we can consider geodesic mappings as trivial almost geodesic
mappings, and we will omit them in further considerations;
they were 
analysed in \cc BerMi=.

Recall that a pseudo-Riemannian space $(M,g)$ is called a
\emph{Ricci-symmetric space}\footnote{In analogy to symmetric spaces that are characterized by parallel
Riemannian curvature tensor} when the Ricci tensor is parallel with respect
to the corresponding Levi-Civita connection ${\nabla}$ of the metric,
${\nabla}{\mbox{\rm Ric}}=0$.
It was proven by Sinyukov \cc Si79=, that the basic partial
differential equations (PDE's) of
${\tilde{\pi}}_1$-mappings of a manifold $(M,\nabla)$ onto
Ricci-symmetric pseudo-Riemannian manifolds $(\bar M,\bar g)$
(of arbitrary signature) can be transformed into
(an equivalent) closed system of PDE's of first order of the Cauchy type.
Hence the solution
(if it exists) depends on a finite set of parameters.
Consequently, for a manifold with a symmetric connection
admitting ${\tilde{\pi}}_1$-mappings onto Ricci-symmetric spaces,
the set of all Ricci-symmetric spaces
$(\bar M,\bar g)$ which can serve as images of the given
manifold $(M,\nabla)$ under ${\tilde{\pi}}_1$-mappings is finite.
The cardinality $r$ of such a set is bounded by the number of
free parameters.

On the other hand, geodesic mappings form a subclass among
${\tilde{\pi}}_1$-mappings (they obey the definition).
Basic equations describing geodesic mappings of manifolds with linear connection
do not form a closed system of Cauchy type (the general solution
depends on $n$ arbitrary functions;
if the given manifold admits geodesic mappings, the cardinality of the
set of possible images is big).
It follows that the conditions (\ref{2.9}) describing
${\tilde{\pi}}_1$-mappings of manifolds, in general,
cannot be transformed into a closed system of Cauchy type.
But if we choose a suitable
subclass of images and restrict ourselves
(for the given manifold) only onto mappings
with co-domain in the apropriate
subclass we might succeed
to get an equivalent closed system of Cauchy type.
If this is the case then the given manifold admits either non (if the
system is non-integrable) or a finite number of
${\tilde{\pi}}_1$-images in the given class.

Our aim is to analyse ${\tilde{\pi}}_1$-mappings of manifolds onto
manifolds with li\-ne\-ar connection in general, and to use the reached results
for examining ${\tilde{\pi}}_1$-mappings
of manifolds onto (pseudo-)Riemannian spaces (in general, without
any restrictive conditions onto the Ricci tensor),
which will generalize the above result by Sinyukov.
In the rest, we will omit ``pseudo''.

All ${\tilde{\pi}}_1$-mappings $f$: $M\to M$ can be described by the following system of
differential equations
\cc Si79, Si82=:
\begin{equation}\label{E*}
\begin{array}{c}
3({\nabla}_Z P(X,Y)+P(Z,P(X,Y))) = \\ \ \\ \ds
\mathop{\sum}_{CS(X,Y)} (R(Y,Z)X-\bar R(Y,Z)X)
+\mathop{\sum}_{CS(X,Y,Z)}a(X,Y)Z.
\end{array}
\end{equation}
In what follows, 
we prefer to express our equalities in local coordinates
(with respect to a map $(U,\varphi)$ on $M$)
since the invariant formulas are rather complicated.
The above formula has the local expression
\begin{equation}
3(P^h_{ij,k}+P^h_{k\a}P^\a_{ij}) = R^h_{(ij)k} - \bar R{}^h_{(ij)k}
+a_{(ij}\d^h_{k)},
\label{4.4}
\end{equation}
where  $P^h_{ij}$, $a_{ij}$, $R^h_{ijk}$, $\bar R{}^h_{ijk}$ are local components
of tensors $P$, $a$, $R$, and $\bar R$.

\section{Properties of the fundamental equations of the ca\-no\-ni\-cal almost geodesic mappings
${\tilde{\pi}}_1$}

Assuming (\ref{4.4}) as a system of PDE's for functions $P^h_{ij}$ on $M$,
the corresponding integrability conditions read
\begin{displaymath}\begin{array}{cc}
\bar R^h_{(ij)[k,\ell]} = R^h_{(ij)[k,\ell]} +
\d^h_{(i} a_{jk),\ell} - \d^h_{(i} a_{j\ell),k} +
3(P^\a_{ij}\bar R^h_{\a k\ell} - P^h_{\a (j}R^\a_{i)k\ell}) -  \\[2mm]
P^h_{\a k}(R^\a_{(ij)\ell} - \bar  R^\a_{(ij)\ell}  \d^\a_{(i} a_{j\ell)}) +
P^h_{\a \ell}(R^\a_{(ij)k} - \bar  R^\a_{(ij)k}  \d^\a_{(i} a_{jk)})
\,.
\end{array}
\end{displaymath}
Passing from $\nabla \bar{R}$ to $\bar{\nabla}\bar{R}$ on the left hand
side we get integrability conditions of the system (\ref{4.4}) in the form
\begin{equation}
\bar R^h_{(ij)[k;\ell]} =
\d^h_{(i} a_{jk),\ell} - \d^h_{(i} a_{j\ell),k} + \Theta^h_{ijk\ell}\, ;
\label{4.6}
\end{equation}
here we denoted
\begin{displaymath}\begin{array}{cc}
\Theta^h_{ijk\ell} = R^h_{(ij)[k,\ell]}  +
3(P^\a_{ij}\bar R^h_{\a k\ell} - P^h_{\a (j}R^\a_{i)k\ell}) -
P^h_{\a k}(R^\a_{(ij)\ell} +  \d^\a_{(i} a_{j\ell)}) + \\[2mm]
P^h_{\a \ell}(R^\a_{(ij)k} +  \d^\a_{(i} a_{jk)}) -
P^\a_{\ell(i} \bar R^h_{|\a |j)k} -
P^\a_{\ell(i} \bar R^h_{j)\a k} +
P^\a_{k(i} \bar R^h_{|\a |j)\ell} +
P^\a_{k(i} \bar R^h_{j)\a \ell}
\end{array}
\end{displaymath}
where ``;"  denotes covariant derivative with respect to $\bar{\nabla}$.

If we apply covariant differentiation with respect to $\bar{\nabla}$
to the integrability conditions (\ref{4.6}) of the system (\ref{4.4}),
and then pass from covariant derivation $\bar{\nabla}$ to ${\nabla}$, we get
\begin{equation}
{\bar R}^h_{(ij)k;\ell m}
-{\bar R}^h_{(ij)\ell;mk}=
{\delta}^h_{(i} a_{jk),\ell m}
-{\delta}^h_{(i} a_{j\ell),km}+T^h_{ijk\ell m} \,,
\label{4.13}
\end{equation}
\smash{where we denoted}
\begin{displaymath}\begin{matrix}
T^h_{ijk\ell m}
={\bar R}^h_{\a mk} {\bar R}^\a_{(ij)\ell} -
{\bar R}^\a_{\ell mk} {\bar R}^h_{(ij)\a} -
{\bar R}^\a_{jmk} {\bar R}^h_{(i\a)\ell} -
{\bar R}^\a_{imk} {\bar R}^h_{(j\a)\ell} -\\
P^h_{m\a}\d^\a_{(i}a_{jk),\ell} -
P^\a_{mj}\d^h_{(i}a_{\a k),\ell} -
P^\a_{mi}\d^h_{(\a}a_{jk),\ell} -
P^\a_{mk}\d^h_{(\a}a_{ij),\ell} -
P^\a_{ml}\d^h_{(i}a_{jk),\a} - \\
P^h_{m\a}\d^\a_{(i}a_{j\ell),k} +
P^\a_{mi}\d^h_{(\a}a_{j\ell),k}+
P^\a_{mj}\d^h_{(i}a_{\a \ell),k} +
P^\a_{mk}\d^h_{(i}a_{j\ell),\a} -
P^\a_{ml}\d^h_{(i}a_{j\a),k} - \\
\Theta^h_{ijk\ell,m}+P^h_{\a m}\Theta^\a_{ijk\ell} -
P^\a_{mi}\Theta^h_{\a jk\ell} - P^\a_{mj}\Theta^h_{i\a k\ell}- P^\a_{mk}\Theta^h_{ij\a \ell}
- P^\a_{m\ell}\Theta^h_{ijk\a } \,.
\end{matrix}
\end{displaymath}
Alternating (\ref{4.13}) in $\ell, m$ we get 
\begin{equation}
\begin{matrix}
\bar R{}^h_{(ij)m;\ell k} - \bar R{}^h_{(ij)\ell;m k} =
\d^h_{(i} a_{jm),k\ell} - \d^h_{(i} a_{j\ell),km}  +
T^h_{ijk[lm]} +
\\
\bar R{}^h_{(i|\a k|} \bar R{}^\a_{j)m\ell} +
\bar R{}^h_{(ij)\a} \bar R{}^\a_{km\ell} -
\bar R{}^\a_{(ij)k} \bar R{}^h_{\a m\ell} +
\bar R{}^h_{\a (i|k|} \bar R{}^\a_{j)m\ell} +
\\
\d^h_{(\a } a_{jk)} R^\a_{i\ell m} +
\d^h_{(\a } a_{ik)} R^\a_{j\ell m} +
\d^h_{(i} a_{j\a )} R^\a_{k\ell m} -
\d^h_{(i} a_{jk )} R^\a_{\a\ell m} \,.
\end{matrix}
\label{4.15}
\end{equation}
Using properties of the Riemannian tensor,
we rewrite (\ref{4.15}) as
\begin{equation}
\bar R{}^h_{im\ell;jk} + \bar R{}^h_{jm\ell;ik}= \d^h_{(i}a_{j\ell),km} -
\d^h_{(i}a_{jm),k\ell}-N^h_{ijk\ell m} \,,
\label{4.16}
\end{equation}
\smash{where the last term is}
\begin{displaymath}
\begin{matrix}
N^h_{ijk\ell m}=T^h_{ijk[\ell m]} +
\bar R{}^\a_{im\ell}  \bar R{}^h_{(\a j) k} +
\bar R{}^\a_{jm\ell}  \bar R{}^h_{(\a i) k} +
\bar R{}^\a_{km\ell}  \bar R{}^h_{(ij) \a} -\\
\bar R{}^h_{\a m\ell}  \bar R{}^\a_{(ij)k} +
\d^h_{(\a} a_{jk)} R^\a_{i\ell m} +
\d^h_{(\a} a_{ik)} R^\a_{j\ell m} +
\d^h_{(\a} a_{ij)} R^\a_{k\ell m} -
a_{(ij} R^h_{k)\ell m} \,.
\end{matrix}
\end{displaymath}
Alternating (\ref{4.16}) over $j,k$ we get
\begin{equation}
\begin{matrix}
\bar R{}^h_{jm\ell;ik} - \bar R{}^h_{km\ell;ij}=
\d^h_{(i}a_{j\ell),km} -
\d^h_{(i}a_{jm),k\ell}-
\d^h_{(i}a_{k\ell),jm} +
\d^h_{(i}a_{km),j\ell} -\\
N^h_{i[jk]\ell m} +
\bar R{}^h_{\a m\ell} \bar R{}^\a_{ikj} +
\bar R{}^h_{i\a \ell} \bar R{}^\a_{mkj} +
\bar R{}^h_{im\a} \bar R{}^\a_{\ell kj} -
\bar R{}^\a_{im\ell} \bar R{}^h_{\a kj} \,.
\end{matrix}
\label{4.18}
\end{equation}
Let us change mutually $i$ and $k$ in
(\ref{4.16}), and then use
(\ref{4.18}).
We evaluate
\begin{equation}
\begin{matrix}
2{\bar{R}}^h_{jm\ell;ik}
=\d^h_{(i} a_{j\ell),km} - \d^h_{(i} a_{jm),k\ell}- \d^h_{(k} a_{jm),i\ell}+ \\
\d^h_{(i} a_{km),j\ell} - \d^h_{(i} a_{k\ell),jm}+\d^h_{(j\ell} a_{k),im}
+{\Omega}^h_{ijk\ell m} \,,
\end{matrix}
\label{4.19}
\end{equation}
where we used the notation
\begin{displaymath}
\begin{matrix}
{\Omega}^h_{ijk\ell m}=
-N^h_{ijk\ell m} + N^h_{k[ij]k\ell m} -
\bar R{}^h_{\a m\ell} \bar R{}^\a_{(kj)i} +
\bar R{}^h_{j\a \ell} \bar R{}^\a_{mik} +
\bar R{}^h_{jm\a} \bar R{}^\a_{\ell ik} -
\\[2mm]
\bar R{}^h_{\a i(j} \bar R{}^\a_{k)m\ell} +
\bar R{}^h_{j\a \ell} \bar R{}^\a_{mik} +
\bar R{}^h_{jm\a} \bar R{}^\a_{\ell ik} -
\bar R{}^h_{\a m\ell} \bar R{}^\a_{ikj} -
\bar R{}^h_{i\a \ell} \bar R{}^\a_{mkj} +
\bar R{}^\a_{im[\ell} \bar R{}^h_{\a] kj} \,.
\end{matrix}
\end{displaymath}
On the left side of (\ref{4.19}),
let us pass from the covariant derivation $\bar{\nabla}$ to $\nabla$:
\begin{equation}
\begin{matrix}
2{\bar R}^h_{jm\ell,ik}=
\d^h_{(i} a_{j\ell),km} - \d^h_{(i} a_{jm),k\ell} - \d^h_{(k} a_{jm),i\ell} \, +  \\
\d^h_{(i} a_{km),j\ell} - \d^h_{(i} a_{k\ell),jm} -\d^h_{(k} a_{j\ell),im}
+S^h_{ijk\ell m}  \,,
\end{matrix}
\label{4.21}
\end{equation}
where
\begin{equation}
\begin{array}{cl}
S^h_{ijk\ell m}=&\Omega^h_{ijk\ell m}- 2\, [
{\bar R}^\a_{jm\ell,i} P^h_{\ell k} -
{\bar R}^h_{\a m\ell,i} P^\a_{jk} -   \\ &
{\bar R}^h_{j\a \ell,i} P^\a_{mk} -
{\bar R}^h_{jm\a ,i} P^\a_{\ell k} -
{\bar R}^h_{jm\ell,\a} P^\a_{ik} + \\ &
({\bar R}^\a_{jm\ell} P^\b_{\a i} -
{\bar R}^h_{\a m\ell} P^\a_{ij} -
{\bar R}^h_{j\a \ell} P^\a_{im} -
{\bar R}^h_{jm\a } P^\a_{i\ell})
P^h_{\b k}  -  \\ &
({\bar R}^\a_{jm\ell} P^h_{\a \b} -
{\bar R}^h_{\a m\ell} P^\a_{\b j} -
{\bar R}^h_{j\a \ell} P^\a_{\b m} -
{\bar R}^h_{jm\a } P^\a_{\b\ell})
P^\b_{i k}  -    \\ &
({\bar R}^\a_{\b m\ell} P^h_{\a i} -
{\bar R}^h_{\a m\ell} P^\a_{\b i} -
{\bar R}^h_{\b\a \ell} P^\a_{i m} -
{\bar R}^h_{\b m\a } P^\a_{i\ell})
P^\b_{j k}  -   \\ &
({\bar R}^\a_{j\b \ell} P^h_{\a i} -
{\bar R}^h_{\a \b\ell} P^\a_{j i} -
{\bar R}^h_{j\a \ell} P^\a_{\b i} -
{\bar R}^h_{j\b \a } P^\a_{i \ell})
P^\b_{k m}  -  \\ &
({\bar R}^\a_{jm\b} P^h_{\a i} -
{\bar R}^h_{\a m\b} P^\a_{j i} -
{\bar R}^h_{j\a \b} P^\a_{m i} -
{\bar R}^h_{jm\a } P^\a_{\b i})
P^\b_{ k\ell}]  \,.
\end{array}
\label{4.22}
\end{equation}

\section{Canonical almost geodesic mappings ${\tilde{\pi}}_1$ onto Riemannian spaces}

Let there exist a ${\tilde{\pi}}_1$-mapping of a manifold $A_n=(M,\nabla)$ onto a Riemannian manifold
${\bar V}_n=(M,\bar g)$ where $\bar g\in T^0_2 M$ is a metric tensor
with components ${\bar g}_{ij}$.
Recall that the Riemannian tensor
$\bar R{}_{hijk}=\bar R{}^\a_{ijk}\bar g{}_{\a h}$
of type $(0,4)$ satisfies
\begin{equation}
{\bar R}_{hijk}+{\bar R}_{ihjk} =0.
\label{5.1}
\end{equation}
In (\ref{4.19}),
let us apply the metric tensor ${\bar g}_{h\beta}$ and then
use symmetrization with respect to $h$ and $j$. According to (\ref{5.1})
we get
\begin{equation}
\begin{matrix}
\bar g_{ih}(a_{m[k,j]l}+a_{l[j,k]m})+
\bar g_{ij}(a_{m[k,h]l}+a_{l[h,k]m})+\\
\bar g_{kh}(a_{m[i,j]l}+a_{l[j,i]})+
\bar g_{kj}(a_{m[i,h]l}+a_{l[h,i]m})+\\
\bar g_{mh}(a_{k[i,j]l}-a_{ij,kl})+
\bar g_{mj}(a_{k[i,h]l}-a_{ih,kl})+
\bar g_{lj}(a_{kh,il}-a_{i(h,k)m})+\\
+2\bar g_{jh}(a_{k(l,i)m}-a_{m(i,k)l})
+\bar g_{lh}(a_{k[j,i]m}-a_{ij,km})
=-\Omega^{\alpha}_{i(j|klm}\bar g_{\a|h)}.
\end{matrix}
\label{5.2}
\end{equation}
Contraction of the last formula  
with the dual tensor ${\bar g}{}^{jh}$ ($\Vert\bar g^{ij}\Vert=\Vert\bar g_{ij}\Vert^{-1}$)
gives
\begin{equation}
\begin{matrix}
a_{kl,im}-a_{im,kl}-a_{km,il}+a_{il,km}
=-\frac{2}{n+1}{\Omega}{}^{\a}_{i\a klm} \,.
\end{matrix}
\label{5.3}
\end{equation}
Let us symmetrize the above formula over $k$ and $l$.
From (\ref{5.3}) we get
\begin{equation}
\begin{matrix}
2a_{kl,im}- 2a_{im,kl}=
2 a_{\a m}R^\a_{lik} + a_{\a i}R^\a_{mlk}+
a_{\a k}R^\a_{mil} + a_{\a l}R^\a_{mik} +\\
\frac2{n+1}(\Omega^\a_{l\a kim} - \Omega^\a_{i\a(kl)m}).
\end{matrix}
\label{5.4}
\end{equation}
Using (\ref{5.3}) and (\ref{5.4}) the equation (\ref{5.2}) reads
\begin{equation}
\begin{matrix}
2\bar g_{ih}(a_{km,jl}-a_{jm,kl})+2\bar g_{ij}(a_{km,hl}-a_{hm,kl})+\\
2\bar g_{kh}(a_{im,jl}-a_{jm,il})+2\bar g_{kj}(a_{im,hl}-a_{hm,il})+\\
\bar g_{mk}(a_{ki,jl}-a_{kj,il}-a_{ij,kl})+
\bar g_{mj}(a_{ki,hl}-a_{kh,il}-a_{ih,kl})+\\
\bar g_{lj}(a_{kh,im}-a_{i(h,k)m})+
\bar g_{lh}(a_{kj,im}-a_{i(k,j)m})=
C_{ijkhl},
\end{matrix}
\label{5.5}
\end{equation}
where
$$ 
\begin{array}{l}
C_{ijkhl}=-\Omega^\a_{i(j|klm} \bar g_{\a|h)}+
\frac2{n+1}  \Omega^\a_{i\a klm} \bar g_{jh}
-\bar g_{kh} a_{\a l}  R^\a_{mij}
+ \\ \quad
\bar g_{ih} (\frac2{n+1}\Omega^\a_{m\a ljk} -
a_{\a k} R^\a_{(ml)j} - a_{\a j} R^\a_{(l|k|m)}-
a_{\a m} R^\a_{lkj} - a_{\a l} R^\a_{mkj}) + \\  \quad
\bar g_{ij} (\frac2{n+1}\Omega^\a_{m\a lhk} -
a_{\a k} R^\a_{(ml)h} - a_{\a h} R^\a_{(l|k|m)} -
a_{\a m} R^\a_{lkh} - a_{\a l} R^\a_{mkh}) + \\  \quad
\bar g_{kh} (\frac2{n+1}\Omega^\a_{m\a lji} -
a_{\a i} R^\a_{(ml)j} - a_{\a j} R^\a_{(l|i|m)} -
a_{\a m} R^\a_{lij} + a_{\a l} R^\a_{mij}) +\\  \quad
\bar g_{kj} (\frac2{n+1}\Omega^\a_{m\a lhi} -
a_{\a i} R^\a_{(ml)h} - a_{\a h} R^\a_{(l|i|m)} -
a_{\a m} R^\a_{lih} + a_{\a l} R^\a_{mih})
.
\end{array}
$$ 
If we contract (\ref{5.5}) with the dual $\bar g{}^{ij}$ of the metric tensor,
use (\ref{5.4}) and the Ricci identity we get
\begin{equation}
\begin{matrix}
a_{km,hl}-a_{kl,hm}=\frac1{2(n+3)}(\bar g_{hm}\mu_{kl} - \bar g_{hl}\mu_{km})+
B_{kmhl},
\end{matrix}
\label{5.7}
\end{equation}
where
$ 
\mu_{km}=a_{\a\b,km} \bar g{}^{\a\b},
$ 
and
$$ 
\begin{array}{l}
B_{kmhl}=C_{\a\b kmhl}\bar g{}^{\a\b} + 3a_{ m\a}R^\a_{lhk} +
\frac32 (a_{h\a }R^\a_{mkl}+ a_{ k\a}R^\a_{mhl}+a_{ l\a}R^\a_{mhk}) +\\
\frac3{n+1}(\Omega^\a_{l\a khm} - \Omega^\a_{h\a (kl)m} ) -
\frac12 (a_{ m\a}R^\a_{lkm}+ a_{k\a }R^\a_{mhl}+a_{ h\a}R^\a_{mkl}+a_{ l\a}R^\a_{mkh}) -\\
\frac1{n+1}(\Omega^\a_{l\a hkm} - \Omega^\a_{k\a (hl)m} ) -
a_{\a (h}R^\a_{k)lm}+
\frac12 (a_{ k\a}R^\a_{lmh} +  a_{ h\a}R^\a_{lkm} +a_{ m\a}R^\a_{lkh}) .
\end{array}
$$ 
Now contract (\ref{5.5}) with $\bar g{}^{ih}$. According to (\ref{5.7}) we get
\begin{equation}
\bar g_{kl}\mu_{jm}- \bar g_{jl}\mu_{km} +\bar g_{km}\mu_{jl}- \bar g_{jkm}\mu_{kl} =
\frac{n+3}{n+1}\, C_{kljm},
\label{5.10}
\end{equation}
where
$$ 
C_{kljm} = C_{\a jkl(m|\b|l)} \bar g{}^{\a\b}
-2(n+1) (B_{k(ml)j}-a_{\a(l}R^\a_{m)jk}+a_{j\a}R^\a_{(m|k|l)} +
 a_{k\a}R^\a_{(lm)j}).
$$ 
Contracting (\ref{5.10}) with $\bar g{}^{k\ell}$ and using the notation
$ 
K=\mu{}_{\a\b}\bar g{}^{\a\b}
$ 
we obtain components of the tensor $\mu$:
\begin{equation}
\mu_{jm}=\frac1n K \bar g_{jm} + \frac{n+3}{n(n+1)} C_{\a\b jm}\bar g{}^{\a\b}.
\label{5.12}
\end{equation}
Using (\ref{5.12}) we can rewrite (\ref{5.7}) in the form
\begin{equation}
a_{km,hl} - a_{hm,kl}=\frac{K}{2n(n+3)} \,
(\bar g_{mh} \bar g_{kl} - \bar g_{lh} \bar g_{km})+ A_{kmhl},
\label{5.13}
\end{equation}
where
$$ 
A_{kmhl} = B_{kmhl} +\frac{1}{2n(n+1)} \,
(\bar g_{hm} C_{\a\b kl}\bar g{}^{\a\b} - \bar g_{hl} C_{\a\b km}\bar g{}^{\a\b}).
$$ 
Combining (\ref{5.5}) and (\ref{5.13}) we get
\begin{equation}
\begin{matrix}
\bar g_{jl}a_{ih,km} + \bar g_{hl}a_{ij,km}-
\bar g_{jm}a_{ih,kl} - \bar g_{hm}a_{ij,kl}= \\
-\frac{K}{n(n+3)} \,
(\bar g_{ih} \bar g_{kl} \bar g_{jm} - \bar g_{ih} \bar g_{km} \bar g_{jl} +
\bar g_{ij} \bar g_{kl} \bar g_{hm} -\bar g_{ij} \bar g_{km} \bar g_{hl} + \\
3\bar g_{kh} \bar g_{il} \bar g_{jm} - 3\bar g_{kh} \bar g_{jl} \bar g_{im} +
3\bar g_{kj} \bar g_{il} \bar g_{hm} - 3\bar g_{lh} \bar g_{jk} \bar g_{im}) + A_{ijkmhl},
\end{matrix}
\label{5.15}
\end{equation}
where we have denoted
$$ 
A_{ijkmhl} = C_{ijkmhl} - 2(\bar g_{i(h}A_{|km|j)l} + \bar g_{k(h}A_{|im|jl)}
- \bar g_{m(h}A_{|ki|j)l} - \bar g_{l(h}A_{|k|j)im)}) .
$$ 
Finally, symetrization of (\ref{5.15}) over the indices $i,j$, followed by contraction with $\bar g{}^{\ell h}$
anables us to express second covariant derivatives of the tensor $a$,
\begin{equation}
a_{ij,km} = \frac K{n(n+3)}\, (\bar g_{ij} \bar g_{km} + 3\bar g_{k(j} \bar g_{i)m})+
A_{(ij)km\a\b}\bar g{}^{\a\b}.
\label{5.17}
\end{equation}
Now we can consider (\ref{5.17}) as the first order system of PDE's of Cauchy type relative to
the tensor $\nabla a$ (i.e.~in $a_{ij,k}$), find the integrability conditions and contract them with
$\bar g{}^{ij}$ and $\bar g{}^{km}$, respectively.
We calculate $\nabla K$,
\begin{equation}
K_{,\b}=\frac{n(n+3)}{n^2+5n-6} \, A_\b ,
\label{5.18}
\end{equation}
where we denoted
$$ 
\begin{matrix}
A_\r=\left[
a_{\a(j,|k} R^\a_{i)m\r} + a_{ij,\a} R^\a_{km\r}-
\frac K{n(n+3)} \,
(\bar g_{ij,[\r} \bar g_{m]k} + \bar g_{ij} \bar g_{k[m,\r}+ \right.\\
\quad
3  \bar g_{kj,[\r} \bar g_{m]i} + 3  \bar g_{kj} \bar g_{i[m,\r]} +
3  \bar g_{ki,[\r} \bar g_{m]j} + 3  \bar g_{ki} \bar g_{j[m,\r]}) + \\
\left.
A_{(ij)k[m|\a\b|,\r]} \bar g{}^{\a\b} +
A_{(ij)k[m|\a\b|}^{\ } \bar g{}^{\a\b}_{\ \ ,\r]}
\right]\,\bar g{}^{ij}\bar g{}^{km}.
\end{matrix}
$$ 
We use $\bar\Gamma{}^h_{ij}=\Gamma{}^h_{ij}+P{}^h_{ij}$ and get
\begin{equation}
\bar g_{ij,k} = P^\a_{ik} \bar g_{\a j} + P^\a_{jk} \bar g_{\a i} .
\label{5.20}
\end{equation}
Assume the tensors $\nabla a$ and $\nabla \bar R$, and denote their components by
$a_{ijk}:=a_{ij,k}$
and
$\bar R{}^h_{ijk\ell}:=\bar R{}^h_{ijk,\ell}$, respectively.
Then (\ref{4.21}) and (\ref{5.17}) take the form
\begin{equation}
\begin{matrix}
2 R^h_{jmli,k} = \d^h_{(i} a_{jl)k,m} -\d^h_{(i} a_{jm)k,l} +
\d^h_{(k} a_{jl)i,m}-\d^h_{(k} a_{jm)i,l} + \\
\d^h_{(i} a_{km)j,l} -\d^h_{(i} a_{kl)j,m} +S^h_{ijklm},
\end{matrix}
\label{5.23}
\end{equation}
\begin{equation}
a_{ijk,m}= \frac K{n(n+3)} (\bar g_{ij} \bar g_{km} + 3\bar g_{k(j} \bar g_{i)m})
+ A_{(ij)km\a\b} \bar g{}^{\a\b},
\label{5.24}
\end{equation}
where
covariant derivatives of the tensor $a{}_{ijk}$ in (\ref{5.23}) are supposed to be expressed according to
(\ref{5.24}), the tensor $S$ was introduced componentwise in~\eqref{4.22}.

The formulas \eqref{4.4}, (\ref{5.18})--(\ref{5.24}) represent a closed system of Cauchy type
for unknown functions
\begin{equation}\label{5.25}
\bar g_{ij}(x),\ P^h_{ij}(x),\ a_{ij}(x),\ a_{ijk}(x),\ K(x),\ \bar R{}^h_{ijk}(x),\  R{}^h_{ijkl}(x),\
\end{equation}
which, moreover, must satisfy a finite set of algebraic conditions
\begin{equation}
\bar g_{[ij]}=P^h_{[ij]}=a_{[ij]}=a_{[ij]k}= \bar R{}^h_{i(jk)}=R{}^h_{i(jk)l}=0, \
det\Vert\bar g_{ij}(x)\Vert\ne 0.
\label{5.26}
\end{equation}
So we have proven:

\begin{theorem}\label{Theo} 
The given manifold $A_n=(M,\nabla)$ admits ${\tilde\pi}_1$-mappings
(i.e. canonical almost geodesic mappings
of type ${\pi}_1$)
onto Riemannian spaces \vnn $=(M,\bar g)$
if and only if there exists solution of the mixed system of Cauchy type
\eqref{4.4}, (\ref{5.18})-(\ref{5.24}), (\ref{5.26}) for the functions (\ref{5.25}).
\end{theorem}

As a consequence of the additional algebraic conditions, we get an upper boundary
for the number $r$ of possible solutions:

\begin{coro}
The family of all Riemannian manifolds ${\bar V}_n$ which
can serve as images of the given manifold $A_n=(M,\nabla)$, depends on at most
$$\frac12n^2(n^2-1)+n(n+1)^2+1$$
parameters.
\end{coro}

The above Theorem
generalizes the result of Sinyukov \cc Si82= already mentioned
as well as his results on geodesic mappings of Riemannian spaces.

\section{Ricci-symmetric and genera\-lized Ricci-sym\-met\-ric spaces}

Under a \emph{Ricci-symmetric manifold}
we mean a manifold $(M,\nabla)$ with linear connection
for which the Ricci tensor is parallel (=\,covariantly constant),
\begin{equation*}
\nabla\mbox{\rm Ric}=0;
\end{equation*}
Ricci symmetric spaces form a particular subclass.

It was proven in \cc Si79=
that the family of all
${\tilde{\pi}}_1$-mappings of a manifold $(M,\nabla)$ onto
Ricci-symmetric ($\bar{\nabla}\bar{\mbox{\rm Ric}}=0$) (pseudo-)Riemannian spaces $(\bar M,\bar g)$
is given by the integrable system
(of Cauchy type)
of partial differentiable equations (in covariant derivatives).
Consequently, given a manifold with a symmetric connection,
the family of all Ricci-symmetric Riemannian spaces $(\bar M,\bar g)$
which can serve as images of the given manifold $(M,\nabla)$
under some ${\tilde{\pi}}_1$-mapping, depends on a finite set of
parameters.

On the other hand, the geodesic mappings form a subset in the set
of ${\tilde{\pi}}_1$-mappings; they obey the definition.
But the basic equations describing geodesic mappings of a manifold
with the linear connection do not form an integrable system
of Cauchy type,
since the general solution depends on $n$ arbitrary functions.
It follows that the conditions (\ref{2.9}) describing
${\tilde{\pi}}_1$-mappings
(i.e.~canonical almost geodesic mappings)
of manifolds do not, in general, induce an integrable
system.

In the following, we consider a particular
case when (\ref{2.9}) can be transformed into an integrable system,
generalizing the results of Sinyukov. Namely,
we will investigate ${\tilde{\pi}}_1$-mappings of a manifold
$(M,\nabla)$ onto the so-called ge\-ne\-ra\-lized
Ricci-symmetric manifolds.
\medskip


A manifold $(M,{\nabla})$ will be called
a \emph{generalized Ricci-symmetric manifold}
if its Ricci tensor satisfies
\begin{equation}
\nabla\mbox{\rm Ric\,}(Y,Z;X)+\nabla\mbox{\rm Ric\,}(X,Z;Y)=0,
\label{5.1a}
\end{equation}
that is,
${\nabla}_X \mbox{\rm Ric\,}(Y,Z)=-{\nabla}_Y\mbox{\rm Ric\,}(X,Z)$.
We do not a priori suppose the Ricci tensor be symmetric.
If $\mbox{\rm Ric\,}$ is symmetric and (\ref{5.1a}) holds
then $\mbox{\rm Ric\,}$
is parallel,
${\nabla}\mbox{\rm Ric\,}=0$, and
$(M,{\nabla})$ is a Ricci-symmetric manifold.
Einstein spaces (Riemannian spaces characterized by the property
that the Ricci tensor is proportional to the metric tensor)
satisfy (\ref{5.1a})
since they satisfy ${\nabla}\mbox{\rm Ric\,}=0$, hence are
generalized Ricci-symmetric.
In this sense, the generalized Ricci-symmetric spaces can be considered
as a certain generalization of Einstein spaces.

\section{Almost geodesic mappings $\tilde\pi_1$ onto genera\-lized Ricci-symmetric manifolds}

Given the $n$-dimensional manifolds ${\A}=(M,{\nabla})$ and
$\bar{\A}=(\bar M,\bar{\nabla})$ with the corresponding
curvature tensors $R$ and $\bar R$, respectively,
all connection-preserving
mappings $f\colon M\to \bar M$ can be described by the system
of differential equations (\ref{E*}),
\cc Si79, Si82,Si83=. These formulas have the local expression (\ref{4.4}).
As we have already proved, from (\ref{E*}) it follows (\ref{4.6}).
Using the Bianci identity we can write (\ref{4.6}) in local coordinates as
\begin{equation*}
{\bar R}^h_{i\ell k;j}+{\bar R}^h_{j\ell k;i}
={\delta}^h_{(i}a_{jk),\ell}-
{\delta}^h_{(i}a_{j\ell),k}+{\Theta}^h_{ijk\ell} \,,
\label{5.6a}
\end{equation*}
where ``;"  denotes the covariant derivative with respect to $\bar{\nabla}$.
Contraction in $h$ and $k$ gives the following equality for covariant
derivatives of components of the Ricci tensor $\bar{\mbox{\rm Ric}}$
of $\bar{\nabla}$:
\begin{equation}
{\bar R}_{i\ell;j} + {\bar R}_{j\ell; i}=
(n+1)a_{ij,\ell}- a_{\ell(i,j)}+{\Theta}^\a_{ij\a\ell} \,.
\label{5.7a}
\end{equation}

\smallskip
In the following let us suppose that the manifold
$(\bar M,\bar{\nabla})$ is a generalized Ricci-symmetric space,
that is, (\ref{5.1a}) holds.
In local coordinates, (\ref{5.1a}) reads
\begin{equation*}
{\bar R}_{ij; k} + {\bar R}_{kj; i}=0.
\label{5.8b}
\end{equation*}
Under this assumption, (\ref{5.7}) reads
\begin{equation}
(n+1)a_{ij,\ell}
-a_{\ell i,j}-a_{\ell j,i}
=-{\Theta}^{\alpha}_{ij\alpha \ell} \, .
\label{5.8a}
\end{equation}
Using symmetrization in $\ell, i$ gives
\begin{equation*}
a_{\ell i,j}+a_{\ell j,i}=-\frac{1}{n}{\Theta}^{\alpha}_{(i|\ell\alpha|j)}
+\frac{2}{n}a_{ij,\ell}  \, .
\label{5.10a}
\end{equation*}
Now (\ref{5.8a}) reads
\begin{equation}
\frac{n^2+n-2}{n}\,
a_{ij,\ell}
=-\,{\Theta}^{\alpha}_{ij\alpha\ell}-\frac{1}{n}
{\Theta}^{\alpha}_{(i|\ell\alpha|j)}.
\label{5.11a}
\end{equation}
Applying the covariant differentiation with respect to $\bar{\nabla}$
to the integrability conditions (\ref{5.4}), followed by passing from
the covariant derivative $\bar{\nabla}$ to ${\nabla}$
on the right hand side, we get
\begin{equation}
{\bar R}^h_{(ij)k;\ell m}
-{\bar R}^h_{(ij)\ell;mk}=
{\delta}^h_{(i} a_{jk),\ell m}
-{\delta}^h_{(i} a_{j\ell),km}+T^h_{ijk\ell m} \,,
\label{5.12a}
\end{equation}
\smash{where}
\begin{equation*}  
\begin{matrix}
T^h_{ijk\ell m}
={\bar R}^h_{\a mk} {\bar R}^\a_{(ij)\ell} -
{\bar R}^\a_{\ell mk} {\bar R}^h_{(ij)\a} -
{\bar R}^\a_{jmk} {\bar R}^h_{(i\a)\ell} -
{\bar R}^\a_{imk} {\bar R}^h_{(j\a)\ell} -\\
P^h_{m\a}\d^\a_{(i}a_{jk),\ell} -
P^\a_{mj}\d^h_{(i}a_{\a k),\ell} -
P^\a_{mi}\d^h_{(\a}a_{jk),\ell} -
P^\a_{mk}\d^h_{(\a}a_{ij),\ell} -
P^\a_{ml}\d^h_{(i}a_{jk),\a} - \\
P^h_{m\a}\d^\a_{(i}a_{j\ell),k} +
P^\a_{mi}\d^h_{(\a}a_{j\ell),k}+
P^\a_{mj}\d^h_{(i}a_{\a \ell),k} +
P^\a_{mk}\d^h_{(i}a_{j\ell),\a} -
P^\a_{ml}\d^h_{(i}a_{j\a),k} - \\
\theta^h_{ijk\ell,m}+P^h_{\a m}\theta^\a_{ijk\ell} -
P^\a_{mi}\theta^h_{\a jk\ell} - P^\a_{mj}\theta^h_{i\a k\ell}- P^\a_{mk}\theta^h_{ij\a \ell}
- P^\a_{m\ell}\theta^h_{ijk\a } \,.
\end{matrix}
\label{5.13a}
\end{equation*}
Alternating (\ref{5.12a}) over $\ell, m$ we obtain
\begin{equation}
\begin{matrix}
\bar R{}^h_{(ij)m;\ell k} - \bar R{}^h_{(ij)\ell;m k} =
\d^h_{(i} a_{jm),k\ell} - \d^h_{(i} a_{j\ell),km}  +
T^h_{ijk[lm]} +
\\
\bar R{}^h_{(i|\a k|} \bar R{}^\a_{j)m\ell} +
\bar R{}^h_{(ij)\a} \bar R{}^\a_{km\ell} -
\bar R{}^\a_{(ij)k} \bar R{}^h_{\a m\ell} +
\bar R{}^h_{\a (i|k|} \bar R{}^\a_{j)m\ell} +
\\
\d^h_{(\a } a_{jk)} R^\a_{i\ell m} +
\d^h_{(\a } a_{ik)} R^\a_{j\ell m} +
\d^h_{(i} a_{j\a )} R^\a_{k\ell m} -
\d^h_{(i} a_{jk )} R^\a_{\a\ell m} \,.
\end{matrix}
\label{5.14a}
\end{equation}
Due to the properties of the Riemannian tensor,
(\ref{5.14a})
can be written as
\begin{equation}
\bar R{}^h_{im\ell;jk} + \bar R{}^h_{jm\ell;ik}= \d^h_{(i}a_{j\ell),km} -
\d^h_{(i}a_{jm),k\ell}-N^h_{ijk\ell m} \,,
\label{5.15a}
\end{equation}
\smash{where}
\begin{equation*}
\begin{matrix}
N^h_{ijk\ell m}=T^h_{ijk[\ell m]} +
\bar R{}^\a_{im\ell}  \bar R{}^h_{(\a j) k} +
\bar R{}^\a_{jm\ell}  \bar R{}^h_{(\a i) k} +
\bar R{}^\a_{km\ell}  \bar R{}^h_{(ij) \a} -\\
\bar R{}^h_{\a m\ell}  \bar R{}^\a_{(ij)k} +
\d^h_{(\a} a_{jk)} R^\a_{i\ell m} +
\d^h_{(\a} a_{ik)} R^\a_{j\ell m} +
\d^h_{(\a} a_{ij)} R^\a_{k\ell m} -
a_{(ij} R^h_{k)\ell m} \,.
\end{matrix}
\label{5.16a}
\end{equation*}
Let us alternate (\ref{5.15a}) over $j$, $k$. We get
\begin{equation}
\begin{matrix}
\bar R{}^h_{jm\ell;ik} - \bar R{}^h_{km\ell;ij}=
\d^h_{(i}a_{j\ell),km} -
\d^h_{(i}a_{jm),k\ell}-
\d^h_{(i}a_{k\ell),jm} +
\d^h_{(i}a_{km),j\ell} -\\
N^h_{i[jk]\ell m} +
\bar R{}^h_{\a m\ell} \bar R{}^\a_{ikj} +
\bar R{}^h_{i\a \ell} \bar R{}^\a_{mkj} +
\bar R{}^h_{im\a} \bar R{}^\a_{\ell kj} -
\bar R{}^\a_{im\ell} \bar R{}^h_{\a kj} \,.
\end{matrix}
\label{5.17a}
\end{equation}
Let us change mutually $i$ and $k$ in
(\ref{5.15a}), and then use
(\ref{5.17a}).
We evaluate
\begin{equation}
\begin{matrix}
2{\bar{R}}^h_{jm\ell;ik}
=\d^h_{(i} a_{j\ell),km} - \d^h_{(i} a_{jm),k\ell}- \d^h_{(k} a_{jm),i\ell}+ \\
\d^h_{(i} a_{km),j\ell} - \d^h_{(i} a_{k\ell),jm}+\d^h_{(j\ell} a_{k),im}
+{\Omega}^h_{ijk\ell m} \,,
\end{matrix}
\label{5.18a}
\end{equation}
where
\begin{equation*}
\begin{matrix}
{\Omega}^h_{ijk\ell m}=
-N^h_{ijk\ell m} + N^h_{k[ij]k\ell m} -
\bar R{}^h_{\a m\ell} \bar R{}^\a_{(kj)i} +
\bar R{}^h_{j\a \ell} \bar R{}^\a_{mik} +
\bar R{}^h_{jm\a} \bar R{}^\a_{\ell ik} -
\\[2mm]
\bar R{}^h_{\a i(j} \bar R{}^\a_{k)m\ell} +
\bar R{}^h_{j\a \ell} \bar R{}^\a_{mik} +
\bar R{}^h_{jm\a} \bar R{}^\a_{\ell ik} -
\bar R{}^h_{\a m\ell} \bar R{}^\a_{ikj} -
\bar R{}^h_{i\a \ell} \bar R{}^\a_{mkj} +
\bar R{}^\a_{im[\ell} \bar R{}^h_{\a] kj} \,.
\end{matrix}
\label{5.19a}
\end{equation*}
On the left hand side of (\ref{5.18a}), let us pass from the covariant derivative with respect to $\bar{\nabla}$ to
the covariant derivative with respect to $\nabla$:
\begin{equation}
\begin{matrix}
2{\bar R}^h_{jm\ell,ik}=
\d^h_{(i} a_{j\ell),km} - \d^h_{(i} a_{jm),k\ell} - \d^h_{(k} a_{jm),i\ell} \, +  \\
\d^h_{(i} a_{km),j\ell} - \d^h_{(i} a_{k\ell),jm} -\d^h_{(k} a_{j\ell),im}
+S^h_{ijk\ell m}  \,,
\end{matrix}
\label{5.20a}
\end{equation}
where
\begin{equation*}\begin{matrix}
S^h_{ijk\ell m}=\Omega^h_{ijk\ell m}- 2\, [
{\bar R}^\a_{jm\ell,i} P^h_{\ell k} -
{\bar R}^h_{\a m\ell,i} P^\a_{jk} -   \\
{\bar R}^h_{j\a \ell,i} P^\a_{mk} -
{\bar R}^h_{jm\a ,i} P^\a_{\ell k} -
{\bar R}^h_{jm\ell,\a} P^\a_{ik} + \\
({\bar R}^\a_{jm\ell} P^\b_{\a i} -
{\bar R}^h_{\a m\ell} P^\a_{ij} -
{\bar R}^h_{j\a \ell} P^\a_{im} -
{\bar R}^h_{jm\a } P^\a_{i\ell})
P^h_{\b k}  -  \\
({\bar R}^\a_{jm\ell} P^h_{\a \b} -
{\bar R}^h_{\a m\ell} P^\a_{\b j} -
{\bar R}^h_{j\a \ell} P^\a_{\b m} -
{\bar R}^h_{jm\a } P^\a_{\b\ell})
P^\b_{i k}  -    \\
({\bar R}^\a_{\b m\ell} P^h_{\a i} -
{\bar R}^h_{\a m\ell} P^\a_{\b i} -
{\bar R}^h_{\b\a \ell} P^\a_{i m} -
{\bar R}^h_{\b m\a } P^\a_{i\ell})
P^\b_{j k}  -   \\
({\bar R}^\a_{j\b \ell} P^h_{\a i} -
{\bar R}^h_{\a \b\ell} P^\a_{j i} -
{\bar R}^h_{j\a \ell} P^\a_{\b i} -
{\bar R}^h_{j\b \a } P^\a_{i \ell})
P^\b_{k m}  -  \\
({\bar R}^\a_{jm\b} P^h_{\a i} -
{\bar R}^h_{\a m\b} P^\a_{j i} -
{\bar R}^h_{j\a \b} P^\a_{m i} -
{\bar R}^h_{jm\a } P^\a_{\b i})
P^\b_{ k\ell}]  \,.
\end{matrix}
\label{5.21a}
\end{equation*}
Let us introduce a $(1,4)$-tensor field
$R^h_{jm\ell i}={\bar R}^h_{jm\ell , i}$. Then we get
\begin{equation}
{\bar R}^h_{jm\ell , i}=R^h_{jm\ell i}  \,.
\label{5.22a}
\end{equation}
From (\ref{5.20a}), the covariant derivative of the tensor (\ref{5.22a}) satisfies
\begin{equation}
\begin{matrix}
2R^h_{jm\ell i , k}=
\d^h_{(i} a_{j\ell),km} - \d^h_{(i} a_{jm),k\ell} - \d^h_{(k} a_{jm),i\ell} +  \\
\d^h_{(i} a_{km),j\ell} - \d^h_{(i} a_{k\ell),jm} +\d^h_{(k} a_{j\ell),im}
+S^h_{ijk\ell m}  \,,
\end{matrix}
\label{5.23a}
\end{equation}
where we used (\ref{5.11a}).

It can be verified that the equations
(\ref{5.2}),
(\ref{5.11a}),
(\ref{5.22a}) and
(\ref{5.23a})
for the functions $P^h_{ij}(x)$, $a_{ij}(x)$, ${\bar R}^h_{ijk}(x)$ and
${R}^h_{ijkm}(x)$
on $(M,\nabla)$ form an integrable system;
the above functions must satisfy also
additional algebraic conditions
\begin{equation}
\begin{matrix}
P^h_{ij}(x)=P^h_{ji}(x),\qquad  a_{ij}(x)=a_{ji}(x),\qquad
{\bar R}^h_{i(jk)}(x) ={\bar R}^h_{(ijk)}(x)=0, \\
R^h_{i(jk)\ell}(x) =R^h_{(ijk)\ell}(x)=0.
\end{matrix}
\label{5.24a}
\end{equation}
So we have succeeded to prove the following
generalization of the result of Sinyukov
\cc Si82, Si83= (we use the  above notation).

\begin{theorem}\label{T2}
Let $(M,\nabla)$ be a manifold with linear connection  and
$(\bar M,\bar{\nabla})$ a generalized Ricci-symmetric manifold.
There
is a ${\tilde{\pi}}_1$ mapping $f\colon M\to \bar{M}$
(i.e.\ a canonical almost geodesic mapping of type $\pi_1$)
if and only if there exist functions $P^h_{ij}(x)$, $a_{ij}(x)$, ${\bar R}^h_{ijk}(x)$ and
${R}^h_{ijkm}(x)$
which satisfy the equations (\ref{5.2}), (\ref{5.11a}),
(\ref{5.22a}), (\ref{5.23a}), and (\ref{5.24a}).
The system of equations (\ref{5.2}), (\ref{5.11a}),
(\ref{5.22a}) and (\ref{5.23a}) forms
 a Cauchy type system of PDE's in covariant derivatives.
\end{theorem}

As a consequence we obtain

\begin{coro}
The family of all generalized Ricci-symmetric manifolds,
which can serve as an image of the given manifold
$(M,\nabla)$ under some ${\tilde{\pi}}_1$-mapping,
depends on at most
\begin{equation}
\frac16\,n(n+1)(2n^3-4n^2+5n+3)
\end{equation}
parameters.
\end{coro}

\end{document}